\documentclass[12pt]{article}
\newtheorem{theorem}{Theorem}
\newtheorem{example}[theorem]{Example}

\usepackage{algorithm}
\usepackage[noend]{algorithmic}

\usepackage{amsmath,amssymb}
\usepackage[margin=3cm]{geometry}
\newcommand{\ffield}[1][]{\ensuremath{\mathbb{F}_{#1}}}

\newcommand{\Aut}[2][\bfield]{\operatorname{Aut}_{#1}(#2)}

\begin{document}

\title{Factoring Ore polynomials over $\mathbb{F}_q(t)$ is difficult\thanks{Research partially supported by grants MTM2013-41992-P and TIN2013-41990-R from the Ministerio de Econom\'{\i}a y Competitividad of the Spanish Government and from FEDER.}
}

\author{Jos\'e G\'omez-Torrecillas$^\star$, F. J. Lobillo$^\star$, Gabriel Navarro$^\dag$ \\
${}^\star$Department of Algebra and CITIC \\
${}^\dag$Department of Computer Sciences and AI, and CITIC \\
University of Granada, Granada, Spain, E18071 \\
\texttt{gomezj@ugr.es, jlobillo@ugr.es, gnavarro@ugr.es}}

\date{}

\maketitle

\begin{abstract}
We show that the effective factorization of Ore polynomials over $\mathbb{F}_q(t)$ is still an open problem. This is so because the known algorithm \cite{GZ} presents two gaps, and therefore it does not cover all the examples.  We amend one of the gaps, and we discuss what kind of partial factorizations can be then computed by using \cite{GZ}.  
\end{abstract}

\section{Statement of the problem}

Let $\mathbb{F}=\mathbb{F}_q$ be the finite field of $q=p^r$ elements, for some prime number $p$, and consider an automorphism $\sigma:\mathbb{F}(t)\to \mathbb{F}(t)$ of the field of rational functions over  $\mathbb{F}$, and $\delta:\mathbb{F}(t)\to \mathbb{F}(t)$ a $\sigma$-derivation. In \cite[Section 4]{GZ}, the algorithm \texttt{Factorization} for factoring out Ore polynomials in $R=\mathbb{F}(t)[x;\sigma,\delta]$ is proposed. \texttt{Factorization} requires as a preprocessing step the computation of the invariant subfield $K$ of $\mathbb{F}(t)$ under $\sigma$. After the standard reduction to the cases of a pure derivation $\mathbb{F}(t)[x;\delta]$ or a pure automorphism  $\mathbb{F}(t)[x;\sigma]$ (\cite[Section 2.1]{GZ}), the latter is further reduced  to considering automorphisms of simpler form, namely, the so called shift and dilation cases, see \cite[Section 2.2]{GZ}. However, such a reduction does not work properly  for the factorization problem. This gap is amended in Section \ref{step1}. 

Even in the cases where $K$ is successfully computed, \texttt{Factorization} does not work for general polynomials.   Given $f\in R$, \texttt{Factorization} first computes a basis of the so-called eigenring $\mathcal{E}(Rf)=\{ u\in R \text{ with } fu\in Rf \text{ and } \deg u<\deg f\}$ as a $K$-algebra, and then the results of \cite{Ivanyos} are required to apply in order to compute either a proper factorization of $f$ or to certify that $f$ is irreducible. However, \texttt{Factorization} fails to compute a factorization for some polynomials. The simpler case is when $f = gh$, where $g,h \in R$ are coprime similar irreducible polynomials. We discuss this problem in Section \ref{step2}.

\section{The computation of the invariant subfield}\label{step1}

In the pure automorphism case $R = \mathbb{F}_q(t)[x;\sigma]$, where $\sigma$ is determined by the linear fractional transformation $t \mapsto \frac{\sigma_1 t + \sigma_2}{\sigma_3 t + \sigma_4} $,  \cite[Section 2.2]{GZ} proceeds by identifying  \(\mathrm{Aut}_{\mathbb{F}_q}{\mathbb{F}_q(t)}\) with \(\operatorname{PGL}(2,\mathbb{F}_q)\) via the map \(\Phi:s = \left( \begin{smallmatrix} \sigma_1 & \sigma_2 \\ \sigma_3 & \sigma_4 \end{smallmatrix} \right) \mapsto \frac{\sigma_1 t + \sigma_2}{\sigma_3 t + \sigma_4} \). Then the authors consider two cases. Case (2) happens when the Jordan form of \(s\) is \(usu^{-1} = \left( \begin{smallmatrix} \alpha & 1 \\ 0 & \alpha \end{smallmatrix} \right)\) for some \(2 \times 2\) non singular matrix $u$ over \(\mathbb{F}_q\). This case can be reduced  to a shift case as explained there. In case (1) the Jordan form of \(s\) is \(usu^{-1} = \left( \begin{smallmatrix} \alpha^q & 0 \\ 0 & \alpha \end{smallmatrix} \right)\), but if the characteristic polynomial of $s$ is irreducible over $\mathbb{F}_q$, then \(\alpha \in \mathbb{F}_{q^2}\) and \(u\) is a \(2 \times 2\) non singular matrix over \(\mathbb{F}_{q^2}\). Let \(\tau\) be the fractional linear transformation corresponding to \(u\). Observe that, in this case, i.e. the characteristic polynomial of $s$ is irreducible over $\mathbb{F}_q$, \(\tau\) is an automorphism of \(\mathbb{F}_{q^2}(t)\), but it cannot be restricted to an automorphism of \(\mathbb{F}_{q}(t)\). Then the automorphism \(\overline{\sigma} := \tau^{-1} \circ \sigma \circ \tau\) (not $\tau \circ \sigma \circ \tau^{-1}$, as claimed in \cite{GZ}, since $\Phi$ is already a group anti-isomorphism) of $\mathbb{F}_{q^2}(t)$ is of dilation type, but it is not an automorphism of $\mathbb{F}_{q}(t)$. The automorphism $\tau^{-1}$ extends canonically to an isomorphism of rings \(\mathbb{F}_{q^2}(t)[x; \sigma] \cong \mathbb{F}_{q^2}(t)[\overline{x}; \overline{\sigma}]\). Actually, the Ore extension \(\mathbb{F}_q(t)[\overline{x}; \overline{\sigma}]\) requires to be defined that \(\overline{\sigma}\) is an automorphism of \(\mathbb{F}_{q}(t)\), which is not the case, so it does not make sense to look for an isomorphism from \(\mathbb{F}_{q}(t)[x; \sigma]\) to \(\mathbb{F}_{q}(t)[\overline{x};\overline{\sigma}]\). 

In the last paragraph of Section 2.2 in \cite{GZ} it is claimed that to factor an \(f \in \mathbb{F}_q (t)[x; \sigma]\) we may factor the polynomial \(\tau^{-1}(f)\). This ring isomorphism is efficiently computable, and thus they assume that \(\sigma\) is either a shift or dilation over \(\mathbb{F}_q (t)\), where \(q\) is redefined as appropriate. This redefinition consists in changing \(q\) by \(q^2\) if the characteristic polynomial of \(s\) is irreducible. As we have seen before, in this case $\tau^{-1}(f)$ only makes sense in $\mathbb{F}_{q^2}(t)[\overline{x}; \overline{\sigma}]$. So, instead of factoring \(f\) as polynomial in \(\mathbb{F}_{q}(t)[x;\sigma]\), they try to factor out \(f\) as an element in \(\mathbb{F}_{q^2}(t)[x;\sigma]\), and the factorizations of \(f \in \mathbb{F}_{q^2}(t)[x;\sigma]\), following \cite{GZ}, correspond to factorizations of $\tau^{-1}(f)$ in $\mathbb{F}_{q^2}(t)[\overline{x}; \overline{\sigma}]$. However, this reduction does not always provide factorizations of \(f\) viewed in \(\mathbb{F}_{q}(t)[x;\sigma]\) as Example \ref{counterexample} shows. 

\begin{example}\label{counterexample}
Let us apply the procedures in \cite{GZ} to the polynomial
\[
f(x) = x^2 + \frac{t^2 + 1}{t} x + (t^2 + t + 1) \in \mathbb{F}_2(t)[x;\sigma]
\]
where \(\sigma(t) = \frac{t+1}{t}\). The characteristic polynomial of \(\left( \begin{smallmatrix} 1 & 1 \\ 1 & 0 \end{smallmatrix} \right)\) is irreducible, hence we have to view \(\sigma\) in  \(\Aut[{\mathbb{F}_4}]{\mathbb{F}_4(t)}\) and \(f \) in \(\mathbb{F}_4(t)[x; \sigma]\). Now we compute the Jordan form:
\[
\left( \begin{matrix} \alpha^2 & 0 \\ 0 & \alpha \end{matrix} \right) = \left( \begin{matrix} \alpha^2 & 1 \\ \alpha & 1 \end{matrix} \right) \left( \begin{matrix} 1 & 1 \\ 1 & 0 \end{matrix} \right) \left( \begin{matrix} 1 & 1 \\ \alpha & \alpha^2 \end{matrix} \right).
\]
As remarked before, \(\tau^{-1}(t) = \frac{t + 1}{\alpha^2 t + \alpha}\) is extended naturally to an isomorphism \(\mathbb{F}_4(t)[x; \sigma] \cong \mathbb{F}_4(t)[\overline{x}; \overline{\sigma}]\), where \(\overline{\sigma}(t) = \alpha t\). Note that \(\overline{\sigma}(t) \not\in \mathbb{F}_2(t)\), and the Ore extension \(\mathbb{F}_2(t)[\overline{x};\overline{\sigma}]\) makes no sense. So we need to factorize \(\tau^{-1}(f) \in \mathbb{F}_4(t)[\overline{x}; \overline{\sigma}]\). Since
\[
\tau^{-1}(f) = \overline{x}^{2} + \left(\frac{\alpha^2 t^{2} + 1}{\alpha^2
t^{2} + \alpha t + 1}\right) \overline{x} + \frac{t}{\alpha^2 t^{2} + \alpha} = \left( \overline{x} + \frac{\alpha t}{\alpha t + \alpha^2} \right) \left( \overline{x} + \frac{\alpha^2}{\alpha t + \alpha^2} \right),
\]
this leads to the following factorization of \(f\),
\[
f = x^2 + \frac{t^2 + 1}{t} x + (t^2 + t + 1) = \left( x + t + \alpha \right) \left( x + t + \alpha^2 \right).
\]
But this is a factorization of \(f\) as Ore polynomial in \(\mathbb{F}_4(t)[x;\sigma]\), which is the one that should be obtained by Algorithm \texttt{Factorization} in \cite[pp. 132]{GZ}. Nevertheless this factorization cannot lead to  any factorization of \(f\) as Ore polynomial in \(\mathbb{F}_2(t)[x;\sigma]\), since \(f\) is irreducible. In order to check this, a bound of \(f\) can be computed by means of \cite{GLN}. Concretely, 
\[
f^* = x^{6} + x^{3} + \frac{t^{6} + t^{5} + t^{3} + t + 1}{t^{4} + t^{2}}.
\]
Since \(\sigma\) has order \(3\), the element \(s = (\sigma^2 + \sigma + 1)(t) = \frac{t^3 + t + 1}{t^2 + t}\) is invariant under \(\sigma\). We have then \(\mathbb{F}_2(s) \subseteq K \subseteq \mathbb{F}_2(t)\). By \cite[Theorem pp. 197]{vanderWaerden:1949}, it follows that \([\mathbb{F}_2(t):\mathbb{F}_2(s)] = 3\), hence \([K:\mathbb{F}_2(s)] = 1\) and \(K = \mathbb{F}_2(s)\). So the center \(C(\mathbb{F}_2(t)[x;\sigma]) = \mathbb{F}_2(s)[x^3]\) by \cite[Theorem 1.1.22]{Jacobson:1996}. Now, since
\[
\frac{t^{6} + t^{5} + t^{3} + t + 1}{t^{4} + t^{2}} = \left( \frac{t^3 + t + 1}{t^2 + t} \right)^2 + \frac{t^3 + t + 1}{t^2 + t} + 1,
\]
it follows that 
\[
f^* = (x^3)^2 + (x^3) + s^2 + s + 1 \in \mathbb{F}_2(s)[x^3],
\]
which is irreducible. Hence \(f\) is irreducible by \cite{GLN}.
\end{example}

Example \ref{counterexample} explains why Algorithm \texttt{Factorization} in \cite[pp. 132]{GZ} fails if \(\sigma\) is not a dilation nor a shift. Nevertheless, this gap can be amended by providing a general method for describing the invariant subfield of $\mathbb{F}(t) = \mathbb{F}_q(t)$ under $\sigma$. Let $\mu$ denote the order of $\sigma$, and let \(K = \ffield(t)^\sigma\). A description of \(K\) appears in \cite{GutierrezSevilla} for any finite subgroup \(H \leq \mathrm{Aut}_{\ffield}{\ffield(t)}\). Let us apply its results to our setting, i.e. \(H = \{1, \sigma, \dots, \sigma^{\mu-1} \}\). In this case, \cite[Algorithm 1]{GutierrezSevilla} can be written as shown in Algorithm \ref{InvSubfield}.

\begin{algorithm}
\caption{Invariant subfield. \cite{GutierrezSevilla}}
\label{InvSubfield}
\begin{algorithmic}
\REQUIRE{$\sigma \in \mathrm{Aut}_{\ffield}{\ffield(t)}$}
\REQUIRE{$e_0, \dots, e_\mu$ the elementary symmetric functions}
\ENSURE{$s \in \ffield(t)$ such that $\ffield(t)^\sigma = \ffield(s)$}
\FOR{$i = 0, \dots, \mu-1$}
\STATE $h_i \gets \sigma^i(t)$
\ENDFOR
\STATE $i \gets 1$
\REPEAT
\STATE $s \gets e_i(h_0, \dots, h_{\mu-1})$
\STATE $i \gets i+1$
\UNTIL{$s \notin \ffield$}
 \RETURN $s$
\end{algorithmic}
\end{algorithm}

Correctness of Algorithm \ref{InvSubfield} is ensured by \cite[Theorem 15]{GutierrezSevilla}. It remains to find a procedure to write any \(f \in \ffield(t)^\sigma\) as a rational function in \(s\), where \(\ffield(t)^\sigma = \ffield(s)\), i.e. we want to find \(g \in \ffield(t)\) such that \(f = g(s)\). This is the Functional Decomposition Problem for univariate rational functions. Although there is a large literature in this FDP, for our purposes, we may refer the approach in \cite{Dickerson:1989}, where the coefficients of \(g\) are computed solving the appropriate system of linear equations. The procedure is better understood with an example:

\begin{example}
We present in this example how the polynomial \(f^*\) in Example \ref{counterexample} is written as an element in \(C(\ffield[2](t)[x;\sigma]) = \ffield[2](s)[x^3]\), where \(s = \frac{t^3 + t + 1}{t^2 + t}\). Although we have computed \(s\) directly in Example \ref{counterexample}, \(s\) is the output of Algorithm \ref{InvSubfield}. In order to do so, we need to find \(g \in \ffield[2](t)\) such that 
\[
\frac{t^{6} + t^{5} + t^{3} + t + 1}{t^{4} + t^{2}} = g \left( \frac{t^3 + t + 1}{t^2 + t} \right).
\]
Since \(\deg(f) = 6\) and \(\deg(s) = 3\), it follows that \(\deg(g) = 2\), i.e. \(g = \frac{g_0 + g_1 t + g_2 t^2}{g_3 + g_4 t + g_5 t^2}\). Then 
\[
\begin{split}
\frac{t^{6} + t^{5} + t^{3} + t + 1}{t^{4} + t^{2}} &= \frac{g_0 + g_1 \frac{t^3 + t + 1}{t^2 + t} + g_2 \left( \frac{t^3 + t + 1}{t^2 + t} \right)^2}{g_3 + g_4 \frac{t^3 + t + 1}{t^2 + t} + g_5 \left( \frac{t^3 + t + 1}{t^2 + t} \right)^2} \\
&= \frac{g_0 (t^2 + t)^2 + g_1 (t^3 + t + 1)(t^2 + t) + g_2 (t^3 + t + 1)^2}{g_3 (t^2 + t)^2 + g_4 (t^3 + t + 1)(t^2 + t) + g_5 (t^3 + t + 1)^2} \\
&= \frac{g_2 + g_1 t + (g_0 + g_2) t^2 + g_1 t^3 + (g_0 + g_1) t^4 + g_1 t^5 + g_2 t^6}{g_5 + g_4 t + (g_3 + g_5) t^2 + g_4 t^3 + (g_3 + g_4) t^4 + g_4 t^5 + g_5 t^6},
\end{split}
\]
which leads to the following linear equations
\begin{align*}
g_2 &= 1 \\
g_1 &= 1 \\
g_0 + g_2 &= 0 \\
g_0 + g_1 &= 0 \\
g_5 &= 0 \\
g_4 &= 0 \\
g_3 + g_5 &= 1 \\
g_3 + g_4 &= 1, 
\end{align*}
whose solution is \(g_0 = g_1 = g_2 = g_3 = 1, g_4 = g_5 = 0\), i.e. \(g = 1 + t + t^2\). We then get 
\[
\frac{t^{6} + t^{5} + t^{3} + t + 1}{t^{4} + t^{2}} = g(s) = 1 + \frac{t^3 + t + 1}{t^2 + t} + \left( \frac{t^3 + t + 1}{t^2 + t} \right)^2
\]
as pointed out in Example \ref{counterexample}.
\end{example}

\section{Finding a zero divisor in a simple Artinian algebra?}\label{step2}

Let $f \in R = \mathbb{F}(t)[x;\sigma,\delta]$ non constant, and consider its eigenring $\mathcal{E}(Rf)$. Algorithm \texttt{Factorization} from \cite{GZ} is based on the fact that either $\mathcal{E}(Rf)$ has a zero-divisor or it is a division ring (\cite[Theorem 3.6]{GZ}). As a consequence, finding a proper factor of $f$ is equivalent to computing a zero-divisor of $\mathcal{E}(Rf)$. The authors of \cite{GZ} claim that this can be effectively done by using the algorithms of \cite{Ivanyos}. In particular, they claim that, if $\mathcal{E}(Rf)$ is semisimple, then  \cite[Theorem 4.4]{Ivanyos} will return a set of primitive ortogonal idempotents of $\mathcal{E}(Rf)$. The problem is that what \cite[Theorem 4.4]{Ivanyos} gives is a complete set of \textbf{central} idempotents of $\mathcal{E}(Rf)$, so, in particular, if $\mathcal{E}(Rf)$ is a simple Artinian ring, then it will return $\{ 1 \}$, and no zero-divisor is computed. However, $\mathcal{E}(Rf)$ has, in general, zero-divisors (it is isomorphic to a matrix ring over a division ring) and $f$ could be reducible. For instance let $f$ be the least left common multiple of two similar irreducible coprime polynomias $g, h \in R$. Here, similar means that $R/Rg \cong R/Rh$ as left $R$--modules. Then $Rf = Rg \cap Rh$ and $R = Rg + Rh$, which implies that $R/Rf \cong R/Rg \oplus R/Rh$. Since $\mathcal{E}(Rf)$ is nothing but the endomorphism ring of $R/Rf$, we see that, in this case, $\mathcal{E}(Rf)$ is isomorphic to a $2 \times 2$ matrix ring and, thus, the procedure described in \cite[Section 4]{GZ} to compute effectively a factorization of $f$ by using \texttt{Factorization} will say wrongly that $f$ is irreducible. Let us mention that pairs  $g,h \in R$ fulfilling these conditions do exist in abundance. The following is an example of minimal size.

\begin{example}
Let $R = \mathbb{F}_2(t)[x;\sigma]$, where $\sigma(t) = t +1$ (so, it is of shift type). The polynomial $f = x^2 + t^2 + t \in R$ factorizes as $f = (x +t)(x+t+1)$. However, $g=x+t, h= x + t +1$ are obviously irreducible, coprime, and they are similar since they have a common bound $f$. If we apply the method of \cite{GZ}, then we get that $f$ is irreducible. 
\end{example}

It is worth to mention that what can be effectively computed by using \cite{GZ} is a factorization of the form $f = q_1 \cdots q_t$, where each $\mathcal{E}(Rq_i)$ is a simple finite-dimensional algebra. However, the procedure from \cite[Section 4]{GZ} cannot be used to further factorize each $q_i$ nor to certify that it is irreducible. Such a partial factorization can be also obtained, avoiding the costly computation of the eigenrings, as soon as a central bound $f^*$ of $f$ is computed and it is factored out in the center (see \cite{GLN}).


\begin{thebibliography}{99}

\bibitem{GZ}
M. Giesbrecht and Y. Zhang
\newblock Factoring and Decomposing Ore Polynomials over $\mathbb{F}_q(T)$,
\newblock {\em Proceedings of the 2003 International Symposium on Symbolic and Algebraic Computation (ISSAC2003)}, 127 -- 134. ACM, New York, NY, USA, 2003.

\bibitem{Ivanyos}
G. Ivanyos and L. R\'onyai and A. Sz\'ant\'o,
\newblock Decomposition of algebras over $\mathbb{F}_q(X_1,\ldots,X_m)$,           
\newblock {\em Applicable Algebra in Engineering, Communication and Computing} 5(2) (1994), 71 -- 90.

\bibitem{GLN}
J. G\'omez-Torrecillas and F.J. Lobillo and G. Navarro, 
\newblock Computing the bound of an Ore polynomial. Applications to factorization,
\newblock preprint 2013. Available at http://arxiv.org/abs/1307.5529

\bibitem{GutierrezSevilla}
J. Gutierrez and D. Sevilla.
\newblock Building counterexamples to generalizations for rational functions of
  Ritt's decomposition theorem.
\newblock {\em Journal of Algebra}, 303(2) (2006), 655 -- 667.
\newblock Computational Algebra.

\bibitem{Dickerson:1989}
Matthew~Thomas Dickerson.
\newblock {\em The functional decomposition of polynomials}.
\newblock PhD thesis, Department of Computer Science, Cornell University,
  Ithaca, NY, 1989.
\newblock PHD.

\bibitem{Jacobson:1996}
N.~Jacobson.
\newblock \emph{{Finite-dimensional division algebras over fields.}}
\newblock {Berlin: Springer}, 1996.

\bibitem{vanderWaerden:1949}
B.~L. van~der Waerden.
\newblock \emph{Modern Algebra}, volume~I.
\newblock Frederick Ungar Publishing Co., 1949.


\end{thebibliography}
\end{document}